\documentclass{article}
\usepackage[a4paper, portrait, margin=1.1811in]{geometry}
\usepackage[english]{babel}
\usepackage[utf8]{inputenc}
\usepackage[T1]{fontenc}
\usepackage{helvet}
\usepackage{algorithm2e}

\usepackage{}
\usepackage{multicol}
\usepackage{etoolbox}
\usepackage{graphicx}
\usepackage{titlesec}
\usepackage{caption}
\usepackage{booktabs}
\usepackage{amsmath}
\usepackage{xcolor} 
\usepackage{comment}
\usepackage[colorlinks, citecolor=cyan]{hyperref}
\usepackage{caption}
\usepackage{subcaption}
\usepackage{cutwin}
\graphicspath{ {./images/} }
\usepackage{scrextend}
\usepackage{fancyhdr}
\usepackage{graphicx}
\newcounter{lemma}

\newcounter{theorem}

\usepackage{amsmath}

\fancypagestyle{plain}{
	\fancyhf{}

}

\makeatletter
\patchcmd{\@maketitle}{\LARGE \@title}{\fontsize{16}{19.2}\selectfont\@title}{}{}
\makeatother

\usepackage{authblk}

\newsavebox\affbox
\author[1]{\textbf{Stewart Jared}}
\author[2]{\textbf{Tokman Mayya}}
\author[3]{\textbf{Bisetti Fabrizio}}
\author[5]{Dallerit Valentin}
\author[4]{\textbf{Diaz-Ibarra Oscar}}
\affil[1,2,5]{ University of California Merced}
\affil[3]{The University of Texas at Austin}
\affil[4]{ Sandia National Laboratories}

\titlespacing\section{0pt}{12pt plus 4pt minus 2pt}{0pt plus 2pt minus 2pt}
\titlespacing\subsection{12pt}{12pt plus 4pt minus 2pt}{0pt plus 2pt minus 2pt}
\titlespacing\subsubsection{12pt}{12pt plus 4pt minus 2pt}{0pt plus 2pt minus 2pt}

\titleformat{\section}{\normalfont\fontsize{10}{15}\bfseries}{\thesection.}{1em}{}
\titleformat{\subsection}{\normalfont\fontsize{10}{15}\bfseries}{\thesubsection.}{1em}{}
\titleformat{\subsubsection}{\normalfont\fontsize{10}{15}\bfseries}{\thesubsubsection.}{1em}{}

\titleformat{\author}{\normalfont\fontsize{10}{15}\bfseries}{\thesection}{1em}{}

\title{\textbf{\huge Variable time-stepping exponential integrators for chemical reactors with analytical Jacobians}\\
	}
\date{}    

\begin{document}

\pagestyle{headings}	
\newpage
\setcounter{page}{1}
\renewcommand{\thepage}{\arabic{page}}

\captionsetup[figure]{labelfont={bf},labelformat={default},labelsep=period,name={Figure }}	\captionsetup[table]{labelfont={bf},labelformat={default},labelsep=period,name={Table }}
\setlength{\parskip}{0.5em}
	
\maketitle
	
\noindent\rule{15cm}{0.5pt}
	\begin{abstract}
        Computational chemical combustion problems are known to be stiff, and are typically solved with implicit time integration methods. A novel exponential time integrator, EPI3V, is introduced and applied to a spatially homogeneous isobaric reactive mixture.  Three chemical mechanism of increasing complexity are considered, and in two cases the novel method can perform similar if not marginally better to a well-known implementation of a BDF implicit method. In one specific case we see relative performance degradation of the EPI3V to the implicit method. Despite this, the novel exponential method does converge for this case. A performance analysis of the exponential method is provided, demonstrating possible avenues for performance improvement.  
	\end{abstract}
\noindent\rule{15cm}{0.4pt}
\section{Introduction}				
	 Combustion is relevant to energy production, transportation, military technology, and most industrial processes. Furthermore, combustion is central to natural events relevant to ecological systems and climate, such as forest fires.  Because of combustion's ubiquity, the ability to model and predict combustion accurately is critical to many engineering and scientific applications. Due to the physical complexity of combustion, numerical simulations have become essential to its study. For example, simulations are used to design high-efficiency, high-performance engines and to predict ignition behavior and pollutant formation \cite{JetSims,FuelEff}.

	The simulation of chemically reactive systems is challenging due their wide range of spatial and temporal scales.  Furthermore, transport of mass, momentum, and energy are tightly coupled to chemical reactions at the molecular scale. For many problems in the low Mach number regime chemical reaction rates are significantly faster than transport processes.  When systems like these are solved numerically it is a common practice to use temporal integration methods with Strang splitting, which advance chemistry and transport separately.  However, integration of the chemical source terms is difficult because it involves a large number of reactions occurring with widely-ranging reaction rates. In other words, while Strang splitting addresses the global stiffness of the problem, integration of the chemistry is still stiff.  The development of efficient time integrators for the chemical source terms appearing in the transport equations for the species concentrations is a critical task in computational combustion.  In this paper, we investigate whether exponential time integration, which has proven efficient in other fields, offers advantages over more established approaches for  chemical kinetics problems in combustion.

	The stiffness of systems of ordinary differential equations that describe the evolution of reactive species and temperature in a homogeneous (zero-d) reactor makes explicit time integration methods impractical since stability constraints on time-step sizes are too severe.  Instead, implicit schemes are typically used.  Commonly employed methods include backward differentiation formulas (BDF)-based integrators \cite{Bisetti1,Solve1,Solve2, Bisetti3}. These algorithms are typically used in conjunction with a modified Newton solver and Krylov-projection-type iterative methods to solve its embedded linear systems \cite{Solve1,Solve2}.  Whether the performance of such implicit methods is satisfactory or not is often predicated on whether an efficient preconditioner can be constructed to accelerate the linear solves. Also, the functional form of the chemical source terms is complex (these source terms and their Jacobians are most often evaluated by software packages such as TCHEM \cite{TCHEMWeb}, Cantera \cite{Cantera}, or Chemkin\cite{ChemKin}), so it is often quite challenging to construct an effective preconditioner \cite{Precondition, Bisetti1}, particularly one that is general enough to be effective across chemical mechanisms. Recently, exponential methods emerged as an efficient alternative to implicit integrators for problems for which an effective preconditioner is not available \cite{EPIRKMayya,EPIRK, Bisetti2}.   
	
	In this paper, we apply a new time adaptive exponential integrator to the simulation of the temporal evolution of chemically reactive and spatially homogeneous systems, i.e. chemical reactors that are described by a system of ordinary differential equations (ODEs). It is found that the novel exponential time integrator will accurately resolve all three chosen chemical mechanisms. The paper is organized as follows.  Section \ref{sec::model} describes the governing differential equations.  In Section \ref{sec::methods}, we discuss the time integration method.  Section \ref{sec::experiments} presents the results and includes a discussion of the comparative performance of legacy implicit and novel exponential methods.  The last Section \ref{sec::conclusion} outlines the conclusions of our study and future directions.  
\color{black}

\section{Governing equations}\label{sec::model}

	We consider a spatially homogeneous chemically reactive system consisting of an ideal gaseous mixture undergoing chemical reactions at constant pressure. The thermochemical state of the mixture is uniquely identified by the mass fractions of each chemical species and temperature.
	
	Chemical species react with each other according to several reactions.  For reaction $j$, the forward reaction rate constant is given given in Arrhenius form by
	\begin{equation}
		f_{j}=A_j T^{\alpha_j} \text{exp}\Bigl( \frac{-E_j}{R T} \Bigr),\label{ArrhenForwar}
	\end{equation}%
	where $A_j$, and $\alpha_j$ are pre-exponential and exponential constants, R is the universal gas constant, and $E_j$ is the activation energy. If reaction $j$ is reversible, then the backwards reaction constant, $b_{j}$, is not zero and is given directly in Arrhenius form (\ref{ArrhenForwar}) or leveraging the equilibrium constant $K_{j}$:
	\begin{equation}
		b_{j}=f_{j} K_{j}.\label{ArrhenReverse}
	\end{equation} 
	The net rate of reaction, representing the number of times that the reaction occurs in the forward direction per unit time per unit volume, is the difference of forward and backward rate constants multiplied by the molar concentration of the species participating in the reaction raised to their stoichiometric coefficients, i.e.
	\begin{equation}
		\mathcal{R}_j = f_{j}\prod_{i=1}^{\mathbf{K}} \mathcal{\chi}_i^{v'_{ji} }-b_{j}\prod_{i=1}^{\mathbf{K}} \mathcal{\chi}_i^{v''_{ji}},
	\end{equation}
	where $v'_{ji}$ and $v''_{ji}$ are the reactant and product stoichiometric coefficients respectively for species $i$ in reaction $j$, and $\mathbf{K}$ is the number of chemical species. 
	
	Production or loss of a species because of reaction is equal to the difference between the forward and backward stoichiometric coefficients times the rate of reaction. We denote the number of reactions as $\mathcal{N}$, and the total amount of species produced is found by summing over all reactions
	\begin{equation}
		\dot{\omega}_i = \sum^{\mathcal{N}}_{j=1} (v''_{ji}-v'_{ji}) \mathcal{R}_j. \label{Production}
	\end{equation}
	The rate of change in temperature is given by the total volumetric heat production divided by thermal capacity of the mixture per unit volume.  The negative sign in equation (\ref{Meq1}) exists because energy is understood with respect to the chemical bonds, not the gas; if these bonds lose energy, the gas' energy increases.  The change in the species mass fractions with respect to time is the net rate of production normalized to mass fractions. The resulting system of ODEs models evolution of the temperature and species mass fractions:
	\begin{subequations}
		\begin{align}
			&\frac{dT}{dt} = - \frac{1}{\rho c_p} \sum_{k=1}^{\mathcal{N}} \dot{\omega}_kH_k W_k, \label{Meq1}\\
			&\frac{dY_i}{dt} = \frac{1}{\rho} \dot{\omega}_iW_i, \hspace{.5 cm} i=1,...,\mathbf{K}. \label{Meq2}
		\end{align}
	\end{subequations}
	The above system of equations include the gas density, $\rho$, the heat capacity at specific pressure, $c_p$, the molar rate of production, $\dot\omega_k$, the specific enthalpy, $H_k$, and the species' molar mass, $W_k$. The system of equations requires a closure for density $\rho$ according to the equation of state for an ideal gas.

	Equations (\ref{Meq1}-\ref{Meq2}) is a general model for a spatially homogeneous isobaric reactive mixture. To finalize the model a specific list of species and reaction parameters must be provided. These are defined through carefully assembled kinetic mechanisms that describe the species present in the gas solution, the chemical reactions that occur, and the species' thermodynamic properties.  The resulting model determines the evolution of the specific gas.
	
	We will use gases consisting of hydrocarbon fuels (C$_x$H$_y$) and oxygen (O$_2$), with additional species such nitrogen (N$_2$) and argon (Ar) accounting for most of the mass \cite{GasComp,Fuels,Surrogate2}. 
	
	
	In order to model the complexity of chemical reactions, and the numerous branching pathways typical of large hydrocarbon oxidation, a large number of species and reactions are required.  We study three hydrocarbon based kinetic mechanisms: GRI3.0 \cite{GRI}, \emph{n}-butane \cite{NBut}, and \emph{n}-dodecane \cite{NDod}, which model the combustion of methane, butane, and dodecane, respectively. Methane, CH$_4$, while the smallest hydrocarbon we study, is found in a wide variety of fuels and is the main component of natural gas. The methane mechanism is detailed and contains 53 species and 325 reactions. The second mechanism is for butane, C$_4$H$_{10}$, and contains 154 species and has 680 reactions. Butane behaves similarly to more complex practical fuels \cite{NBut} and is a component of gasolines \cite{GasComp}. The last mechanism is \emph{n}-dodecane, which is the largest hydrocarbon, C$_{12}$H$_{26}$, we study.  It is a component of kerosene and some jet fuels \cite{Kerosene}, it contains 105 species and 420 reactions.

\section{Methods \& Implementation}\label{sec::methods}
	For ease of notation, we organize the thermochemical state variables in a vector ordered with temperature followed by species in the same order as they appear in the chemical mechanism:
	\begin{equation}
		y(t) = [T(t), \hspace{.25cm} Y_1(t), \hspace{.25cm}  ... \hspace{.25cm} ,\hspace{.25cm}  Y_{N_s}(t)]^T.\label{State}
	\end{equation}
	Denoting the initial gas state of the mixture $y(t_0)=y_0$, one obtains the initial value problem
	\begin{subequations}
	\begin{align}
		\frac{d y(t)}{dt} &=F(y(t)), \\
		y(t_0)&=y_0\label{MeqFinal}.
	\end{align}
	\end{subequations}
	Time is discretized as $[t_0, t_1, \dots , t_m]$, where $t_{n+1} = t_{i} + h_i$, and $h_i$ is the time step size. Approximations of the state, right-hand-side function and Jacobian at time $t_n$ are denoted as $y_n \approx y(t_n)$, $F_n \approx F(y(t_n))$, $J_n \approx J(y(t_n))$, respectively.

	 We employ TCHEM \cite{TCHEMWeb} in order compute the chemical source terms and Jacobians. Given the mass fractions, temperature, and gas pressure, TCHEM returns both $F_n$ and $J_n$.  The Jacobian can be computed either via finite differences, or analytically.  We use the analytical version to avoid the inaccuracies associated with approximations. Moreover, it is well known that computation of the Jacobian via finite differences is inefficient and computationally expensive for large mechanisms \cite{KimTChem}.
	 
	  The simulation of combustion processes requires the solution of nonlinear and stiff systems of ordinary differential equations, which can be large in size depending on the chemical mechanism. Exponential propagation iterative methods of Runge–Kutta type (EPIRK) have been shown to perform efficiently for a range of large scale stiff systems \cite{KIOPS}, including reaction-diffusion models \cite{EPIRKMayya}. Because of the record of success of EPIRK methods and the requirement that the time step size be adaptive in order to simulate ignition problems, we extend the EPIRK framework to  a novel time-adaptive third-order EPIRK method with an embedded second order scheme for error estimation. The details of the derivation of order conditions and their solution for constructing a particular scheme can be found in \cite{EPIRKMayya}. The same approach is used to formulate the following EPIRK integrator, EPI3V:
		\begin{subequations}
			\begin{align}
				\begin{split}
					Y_1 &= y_n +  \varphi_1\Bigl(\tfrac{3}{4} h_n J_n \Bigr) h_nF_n , \label{Stage}
				\end{split}\\
				\begin{split}
					R(z) &= f(z)-F_n - J_n (z-y_n), \label{Remainder}
				\end{split}\\
				\begin{split}
					y_{n+1}&= y_n + \varphi_1(h_n J_n)  h_nF_n +  \varphi_3(h_n J_n) 2 h_n R(Y1). \label{Step}
				\end{split}
			\end{align}
		\end{subequations}
		The $\varphi$-functions are
		\begin{equation}
			\varphi_1(z) = \frac{e^z-1}{z}, \hspace{.5cm}  \varphi_3(z)=\frac{e^z - \frac{1}{2} z^2 -z - 1}{z^3}.
		\end{equation}
		The method above uses matrix arguments for the $\varphi$-functions; computing approximations of the product of exponential-like matrix functions and vectors of type $\varphi_k(A)v$ is the largest computational expense of exponential integrators. Systems of $N$ ordinary differential equations that model realistic physical processes result in large exponential matrices of size $N \times N$ that make the evaluation of the $\varphi$-functions prohibitively expensive with traditional approximations like  Padé \cite{Pade} or Taylor expansions.	Our EPIRK method was designed to leverage KIOPS, which is an adaptive Krylov-projection algorithm designed to estimate $\varphi$-functions \cite{KIOPS}. In KIOPS, an augmented matrix $\tilde{A}$ is used to express the linear combination of $\varphi$-functions as:
		\begin{equation}
			w(\tau)= \sum^{p}_{j=0}\tau^j\varphi_j(\tau A)b_j=e^{\tau \tilde{A}} v.
		\end{equation}
		 A sub-stepping procedure is then employed to estimate the successive products of matrix exponentials and vectors by iteratively letting $\tau = \tau_1 + \tau_2 + ... + \tau_M$:
		 \begin{equation}
			e^{\tau \tilde{A}} v= e^{(\tau_0 + \tau_1 + ... + \tau_M) \tilde{A}} v = e^{\tau_0 \tilde{A}} e^{\tau_1 \tilde{A}} ... e^{\tau_M \tilde{A}} v.
		\end{equation}		 
		Each product $e^{(\tau_l \tilde{A})}v_{\tau_l}$ is approximated with a Krylov projection in the KIOPS algorithm. With the length of Krylov basis being $m$,  $V$ an $N \times m $ matrix with Krylov basis vectors $v_{\tau_{i}}$ as its columns, and $H$ an $m \times m$ matrix so that $H_{ij}=(\tilde{A}v_{\tau_{r,j}})^Tv_{\tau_{r,i}}$ we have the projection
		\begin{equation}
			e^{\tau_i \tilde{A}} v_{\tau_i}\approx V e^{\tau_i H}V^T v_{\tau_i},
		\end{equation}  
		where $e^{\tau_l H}$ is approximated using Padé with a squaring and scaling algorithm \cite{Pade}.
		
		
		We use the Exponential Propagation Integrators Collection (EPIC) C++ package \cite{EPICSoftware}, which includes implementation of EPIRK methods and KIOPS and allows for easy implementation of new methods.  EPIC provides a linear combination of products of $\varphi$-functions and vectors or a single product estimated at various scalar multiples of the $\varphi$-function's arguments.  In order to obtain linear combination of $\varphi$-function vector products, the user provides a matrix $A$ and vectors $b_i$. EPIC then uses KIOPS to approximate
		\begin{equation}
			\varphi_0 (A) b_0+ \varphi_1 (A) b_1 + ... + \varphi_p(A)b_p. \label{Singlecall}
		\end{equation}
		The user can also provide a single $b_i$ and a set of intermediate time points
		\begin{equation}
			[T_1, \dots, T_M], \hspace{.25cm} T_{j+1}>T_j, \hspace{.25cm} T_j\in (0,1).
		\end{equation}
		The KIOPS algorithm allows the time integrator to stop at each $T_j$ and save the values:
		\begin{equation}
			\varphi_i(T_1 A) b_i, \hspace{.25cm}\varphi_i(T_2 A) b_i, \hspace{.25cm} ... \hspace{.25cm}, \hspace{.25cm}\varphi_i(T_M A) b_i ,\hspace{.25cm}\varphi_i(A) b_i.\label{Multicall}
		\end{equation}
		
		We accomplish all $\varphi$-function approximations with two calls to KIOPS. The first call estimates both
		\begin{equation}
			\varphi_1\Bigl(\tfrac{3}{4} h_n J_n \Bigr) h_nF_n , \text{ and }  \varphi_1(h J_n)h_n F_n.
		\end{equation}
		The second call estimates
		\begin{equation}
			\varphi_3(h_n J_n) 2 h_n R(Y1). \label{LTE}
		\end{equation}
		Two separate calls are necessary in order to obtain the local truncation error in support of the adaptive time step size selection algorithm. The lower order exponential Euler method is:
		\begin{equation}
			y_{n+1}=y_n + h_n \varphi_1(hJ_n)F_n.\label{ExpEuler}
		\end{equation}
		Thus, by subtracting the right hand side of equation (\ref{ExpEuler}) from the right hand side of equation (\ref{Step}), we obtain a local truncation error estimate:
		\begin{equation}
			\varphi_3(h_n J_n) 2 h_n R(Y1). \label{LTE_Est}
		\end{equation}  
		This quantity is obtained by the second call to the KIOPS algorithm (\ref{LTE}).
		
		We implemented a standard adaptive controller from Wanner et al. \cite{Controller} in order to create a time-adaptive method. After a step is calculated, the local truncation error (\ref{LTE}) is compared with the controller's tolerance. If not within tolerance, and the step is rejected, the time step is adjusted and the process is repeated.  Once the tolerance is achieved, the step is accepted, the time step is adjusted, and the iteration proceeds to the next step. We also implement slight modifications due to the specific features of an ignition process.  Step sizes change dramatically during the chain-branching phase of ignition and shortly thereafter. The following constraints limit the change between step sizes: $h_{old}$, $h_{new}$, and the estimated new step size $\hat{h}$
		\begin{equation}
			\begin{cases} 
				h_{new} = 2 \hat{h}, & \hat{h} > 100 h_{old} \\
				h_{new} = \frac{1}{100} \hat{h}, & \hat{h} < 1000 h_{old}. \\
			\end{cases}
		\end{equation}
	\\   
	Implicit methods have shown to be an effective for solving stiff systems of arising from the modeling of homogeneous chemically reacting systems for over thirty years \cite{Solve1,Solve2}. The CVODE package from Lawrence Livermore National Laboratory is widely used for solving general systems of ODEs.  This C++ package implements a variable-coefficient ODE (VODE) solver offering adaptivity of the step size \cite{CVODE}.  The user can select linear and non-linear solvers from a list of available options.  Because the proposed EPI3V method uses Krylov projection methods, we configure CVODE to allow for the most informative comparison possible. We choose the BDF non-linear solver and the SPGMR (Scaled Preconditioned Generalized Minimized Residual) linear solver.  Additionally, we force the Jacobian to be evaluated once each time step and fix the maximum order to 3.
	
	The EPIC package uses NVector data structures from SUNDIALS for vector operations. The current version of EPIC is compatible with SUNDIALS v5.7.

\section{Numerical experiments}\label{sec::experiments}

Our novel time integration method is compared to that implemented in CVODE by simulating the ignition of mixtures of air and hydrocarbons. Ignition is a fundamental process in combustion, whereby the gaseous species in the mixture undergo accelerating exothermic chemical reactions. As the energy in the chemical bonds of the fuel is converted into sensible enthalpy, it contributes to an increase in the mixture temperature, which in turn leads to an acceleration of the rate of reaction. Thus, ignition is characterized by a sudden and abrupt exponential increase of temperature and rate of chemical reactions, which are accompanied by a corresponding depletion of fuel and oxidizer and the formation of products of combustion, i.e. water and carbon dioxide. Once either the fuel or oxidizer are exhausted, the mixture reaches an elevated equilibrium temperature and composition, which no longer vary in time.

In this work, we hasten the onset of ignition by setting the initial temperature at or above 1000 Kelvin, which is sufficient to induce the thermal decomposition of molecular oxygen into its O atoms, which commence ignition by attacking the fuel molecule.

For our numerical experiments we ensure that the parameter values for simulations, including temperature, pressure and initial mass fractions, are sufficient for auto-ignition and lean fuel mixtures.  The lean mixtures mean that stoichiometrically there is more oxygen than fuel which will ensure the combustion terminates when the fuel is consumed, preventing reactions between products of combustion with the fuel hydrocarbons.  Experiments are run through the initial buildup phase into the ignition phase until the steady state is achieved. Table \ref{ICs} shows the temperature, mass fractions and final simulation times for each kinetic mechanism.	

	Data sets are generated by executing simulations of ignition and storing the state vector at the end of the time integration interval. The same simulations are run repeatedly with different absolute and relative tolerances in order to obtain a set of solutions of increasing accuracy for each kinetic mechanism. The error at the final time is computed with respect to a reference solution generated using CVODE with tight relative and absolute tolerances.	\\
		\begin{figure}[!h]
			\centering
			\begin{subfigure}[b]{0.45\textwidth}
				\includegraphics[width=\textwidth]{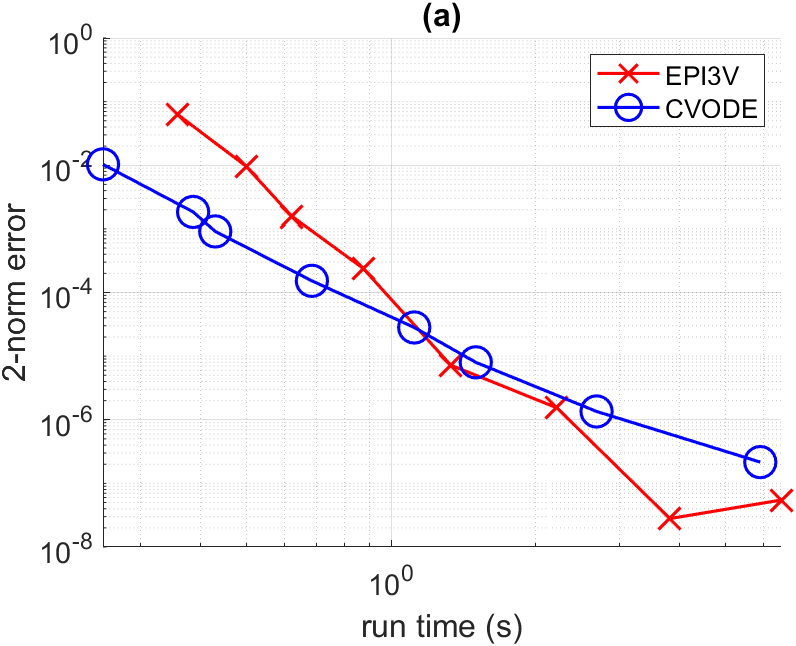}
			\end{subfigure}  
			\begin{subfigure}[b]{0.45\textwidth}
				\centering
				\includegraphics[width=\textwidth]{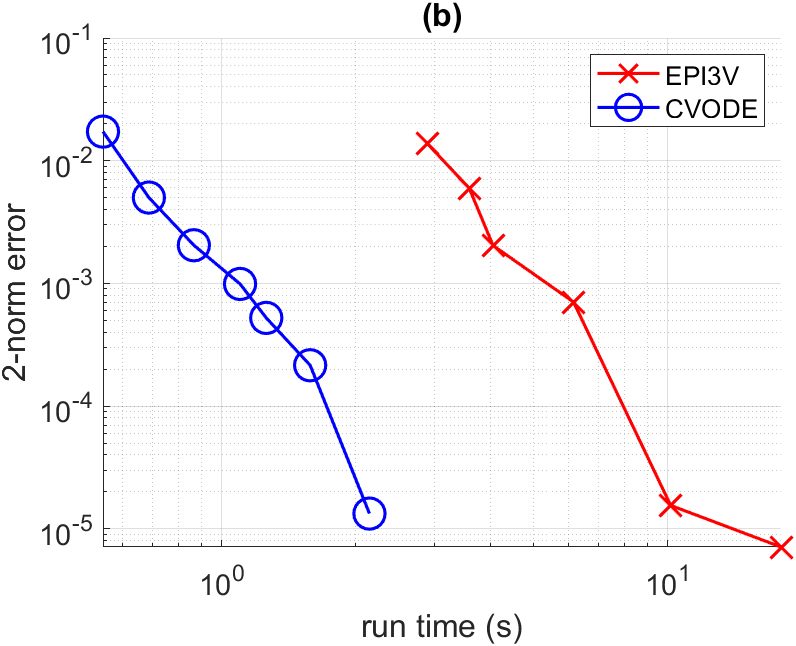}
			\end{subfigure}  
			\begin{subfigure}[b]{0.45\textwidth}
				\includegraphics[width=\textwidth]{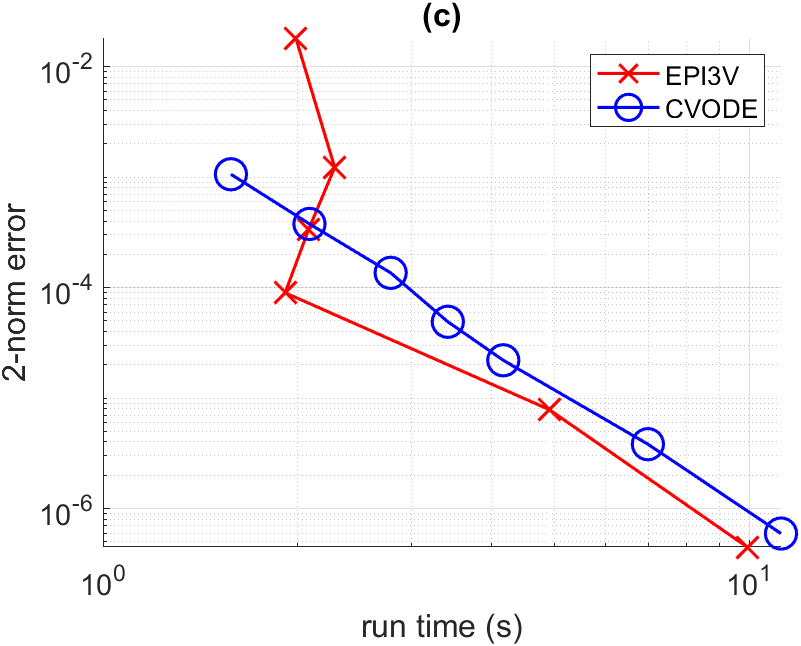}  
			\end{subfigure}
			\caption{Precision diagrams comparing the CPU run time against the 2-norm error of EPI3V versus CVODE. The plots show GRI3.0 (a), $\emph{n}$-dodecane (b), and $\emph{n}$-butane (c) respectively.}
			\label{Precision}
		\end{figure}
		
\begin{table}
	\caption {Experiment configurations.} \label{tab:title} 
	\centering
	\begin{tabular}{ |p{1.5cm} p{3.8cm}|  }
		\hline
		\multicolumn{2}{|c|}{GRI 3.0 Mechanism initial values} \\
		\hline
		Item & Value \\
		\hline
		Kelvin		       &   1000\\
		CH$_4$                &   0.0548\\
		O$_2$                 &   0.2187\\
		Ar                 &   0.0126\\
		N$_2$                 &   0.7137\\
		\hline
	\end{tabular}
	\label{IDTGri}
	\quad
	\begin{tabular}{ |p{1.5cm} p{3.8cm}|  }
		\hline
		\multicolumn{2}{|c|}{$\emph{n}$-butane Mechanism initial values} \\
		\hline
		Item &  Value \\
		\hline
		Kelvin   &   1200\\
		O$_2$                 &     0.2173\\
		C$_4$H$_{10}$              &     0.0607\\
		Ar                 & 	 0.0125\\
		N$_2$                 &  	 0.7092\\
		\hline
	\end{tabular}
	\label{IDTNBut}
	\centering
	\hspace{0.0125cm}
	\begin{tabular}{ |p{1.5cm} p{3.8cm}|}
		\hline
		\multicolumn{2}{|c|}{$\emph{n}$-dodecane Mechanism initial values} \\
		\hline
		Item &  		Value \\
		\hline
		Kelvin   		&   1200\\
		O$_2$                 	&   0.2169\\
		C$_{12}$H$_{26}$             	& 	0.0624\\
		N$_2$                 	&  	0.7080\\
		\hline
	\end{tabular}
	\label{IDTNDod}\centering
	\quad
	\begin{tabular}{ |p{1.8cm} p{3.5cm}|}
		\hline
		\multicolumn{2}{|c|}{Mechanism final times} \\
		\hline
		Mechanism  		&  	Final time (s) \\
		\hline
		GRI 3.0    			&   1.2\\
		$\emph{n}$-butane            	&   $2 \cdot 10^{-3}$	\\
		$\emph{n}$-dodecane           & 	$5 \cdot 10^{-4}$   \\
		\hline
	\end{tabular}
	\label{ICs}
\end{table}
		\begin{table}[h!]
	\caption{Absolute and relative tolerances used to generate precision diagrams in figure \ref{Precision}.}
	\label{ExpTol}
	\centering
	\begin{tabular}{ |p{2.0cm} p{10.0cm}|}
		\hline
		\multicolumn{2}{|c|}{GRI3.0} \\
		\hline
		Method 		&   (Absolute Tolerance, Relative Tolerance)   \\ 
		\hline 
		\rule{0pt}{3ex} EPI3V   			&   ($10^{-10}$, $2 \cdot 10^{-2})$,  ($10^{-10}$, $3 \cdot 10^{-3}$), ($10^{-10}$, $5 \cdot 10^{-4}$), ($10^{-10}$, $10^{-4}$),  \\
		 \rule{0pt}{3ex}&($10^{-10}$, $10^{-5}$), ($10^{-11}$, $5\cdot10^{-6}$), ($10^{-12}$, $2\cdot10^{-6}$), ($10^{-13}$,$8 \cdot 10^{-7}$).\\  \hline \rule{0pt}{3ex}
		CVODE            	&   ($10^{-7}$, $10^{-5}$), ($10^{-8}$, $10^{-6}$), ($10^{-8}$, $10^{-7}$), ($10^{-9}$,$10^{-8}$),\\
		\rule{0pt}{3ex}& ($10^{-10}$, $10^{-10}$), ($10^{-10}$, $10^{-11}$), ($10^{-11}$, $10^{-11}$), ($10^{-11}$, $10^{-12}$). \\
		\hline
	\end{tabular}
	\vspace{.5cm}
	
	\begin{tabular}{ |p{2.0cm} p{10.0cm}|}
		\hline
		\multicolumn{2}{|c|}{$\emph{n}$-butane} \\
		\hline
		Method 		&   (Absolute Tolerance, Relative Tolerance)   \\ 
		\hline
		\rule{0pt}{3ex}EPI3V   			&   ($10^{-7}$, $4 \cdot 10^{-4}$), ($4\cdot 10^{-8}$, $2 \cdot 10^{-4}$), ($2 \cdot 10^{-8}$, $10^{-4}$),   \\
		\rule{0pt}{3ex}&($10^{-8}$, $10^{-4}$), ($10^{-10}$, $10^{-5}$), ($10^{-11}$, $10^{-6}$).\\  \hline
		\rule{0pt}{3ex}CVODE            	&   ($10^{-10}$, $10^{-5}$), ($10^{-10}$, $10^{-6}$), ($10^{-10}$, $10^{-7}$),  \\
		&  ($10^{-10}$,$10^{-9}$), ($10^{-11}$, $10^{-10}$), ($10^{-12}$, $10^{-11}$). \\
		\hline
	\end{tabular}
	\vspace{.5cm}
	
	\begin{tabular}{ |p{2.0cm} p{10.0cm}|}
		\hline
		\multicolumn{2}{|c|}{$\emph{n}$-dodecane} \\
		\hline
		Method 		&   (Absolute Tolerance, Relative Tolerance)   \\ 
		\hline
		\rule{0pt}{3ex}EPI3V   			&   ($10^{-5}$, $5 \cdot 10^{-4}$), ($2 \cdot 10^{-6}$, $2 \cdot 10^{-4}$), ($10^{-7}$, $5 \cdot 10^{-5}$),   \\
		\rule{0pt}{3ex}&($10^{-8}$, $10^{-5}$), ($5\cdot10^{-9}$, $5 \cdot 10^{-6}$), ($10^{-10}$, $10^{-6}$), ($10^{-11}$, $10^{-7}$).\\  \hline
		\rule{0pt}{3ex}CVODE            	&   ($10^{-8}$, $10^{-3}$), ($10^{-8}$, $10^{-4}$), ($10^{-8}$, $10^{-5}$), ($10^{-8}$,$10^{-6}$), \\
		\rule{0pt}{3ex}&  ($10^{-8}$, $10^{-7}$), ($10^{-8}$, $10^{-9}$), ($10^{-8}$, $10^{-10}$), ($10^{-8}$, $10^{-11}$). \\
		\hline
	\end{tabular}
\end{table}

		\begin{table}[h!]
			\caption{Tolerances used to generate the reference solutions in figure \ref{Precision}.}
			\label{RefTol}
			\centering
			\begin{tabular}{ |p{3.0cm} p{3.0cm} p{3.0cm}|}
				\hline
				\multicolumn{3}{|c|}{Reference tolerances} \\
				\hline
				Mechanism name 		&  	Absolute tolerance & Relative tolerance \\ 
				\hline
				\rule[1ex]{0pt}{2ex}GRI 3.0    			&   $10^{-13}$ & $10^{-13}$\\  
				\rule[1ex]{0pt}{1ex}NButane            	&   $10^{-12}$ & $10^{-12}$	\\ 
				\rule[1ex]{0pt}{1ex}NDodecane           & 	$10^{-10}$ & $10^{-10}$\\
				\hline
			\end{tabular}
		\end{table}

		We compare the performance of the EPI3V method against that of CVODE by plotting precision diagrams (CPU time versus a measure of accuracy, here the 2-norm of the error vector) for both methods in figure \ref{Precision}. Tolerances were chosen in order to generate error values for the two methods.  Table \ref{ExpTol} contains the tolerances selected for each numerical experiment.  Table \ref{RefTol} contains the tolerances which generate each experiment's CVODE reference solution.

For the GRI3.0 mechanism, we see a modest advantage in performance for CVODE if loose tolerances are used. However, at tighter tolerances that yield errors below $10^{-5}$ the EPI3V method outperforms CVODE.  With the butane mechanism, we observe a similar relative performance between the two methods. For loose tolerances, the computational time spent integrating the method is insensitive to accuracy, although the wall-clock time is larger than with CVODE until the error is approximately $10^{-5}$. For errors lower than $10^{-5}$, the EPI3V method becomes slightly faster than CVODE. It is important to note that while CVODE is an established code with decades of optimization, our EPI3V implementation is rather new; both software and algorithmic optimizations are ongoing and improvements are expected. For example, significant computational savings were obtained for EPIRK methods recently as the exponential matrix functions evaluations transitioned from straightforward Krylov projection to adaptive Krylov method {\it phipm} \cite{PHIPM} and later to KIOPS \cite{KIOPS}.  In fact, in addition to improvements of the methods' parameters new algorithms may be beneficial for approximating exponential matrix functions for select problems.  Our third test problem, $\emph{n}$-dodecane, illustrates the importance of research in this direction.  

		Unlike the other two mechanisms, simulation of ignition with the $\emph{n}$-dodecane mechanism presented a challenge for the EPI3V method.  Like in the other experiments, we verified the EPI3V method generated the correct solution.  However, unlike for the two previous cases, EPI3V was consistently slower than CVODE by an order of magnitude.  To explore potential causes we considered the spectrum of the Jacobian matrix for all cases.

		\begin{figure}[h!]
			\centering
			\begin{subfigure}[b]{0.45\textwidth}
				\includegraphics[width=\textwidth]{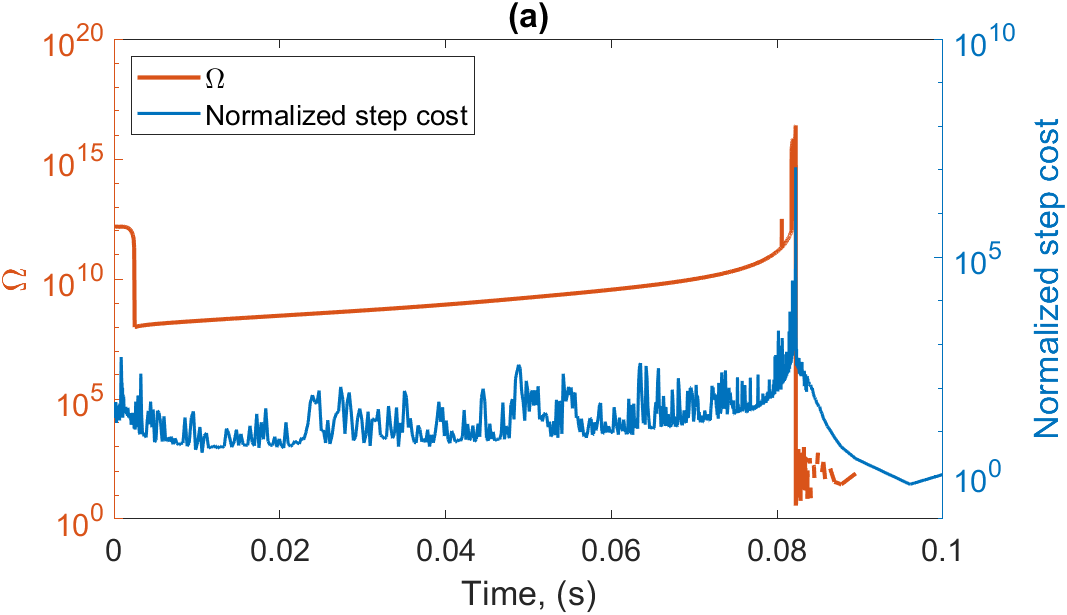}
			\end{subfigure}  
			\begin{subfigure}[b]{0.45\textwidth}
				\includegraphics[width=\textwidth]{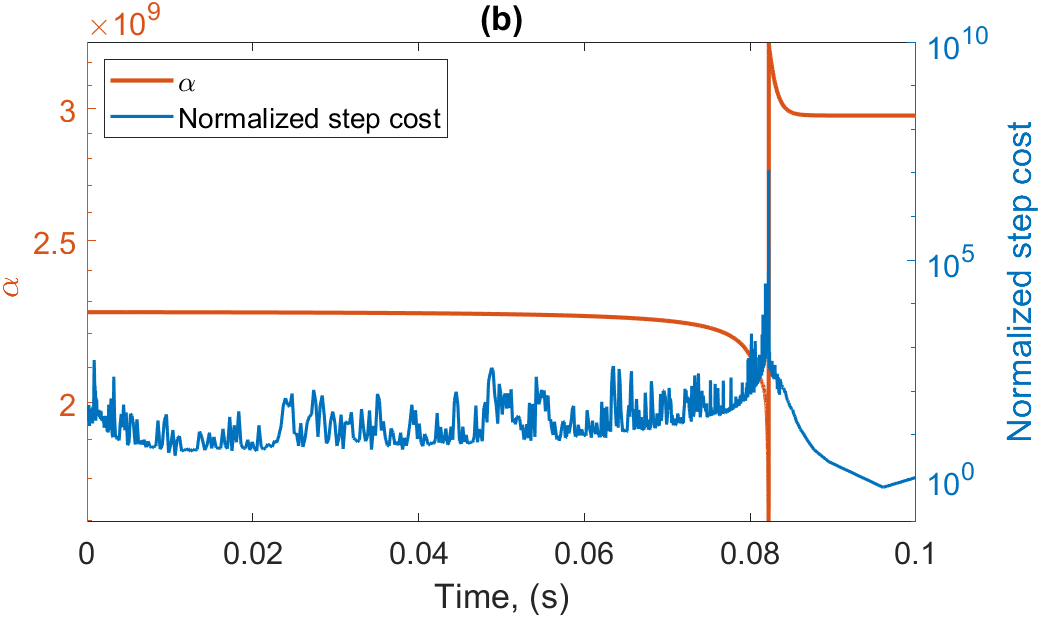}
			\end{subfigure}
			\begin{subfigure}[b]{0.45\textwidth}
				\centering
				\includegraphics[width=\textwidth]{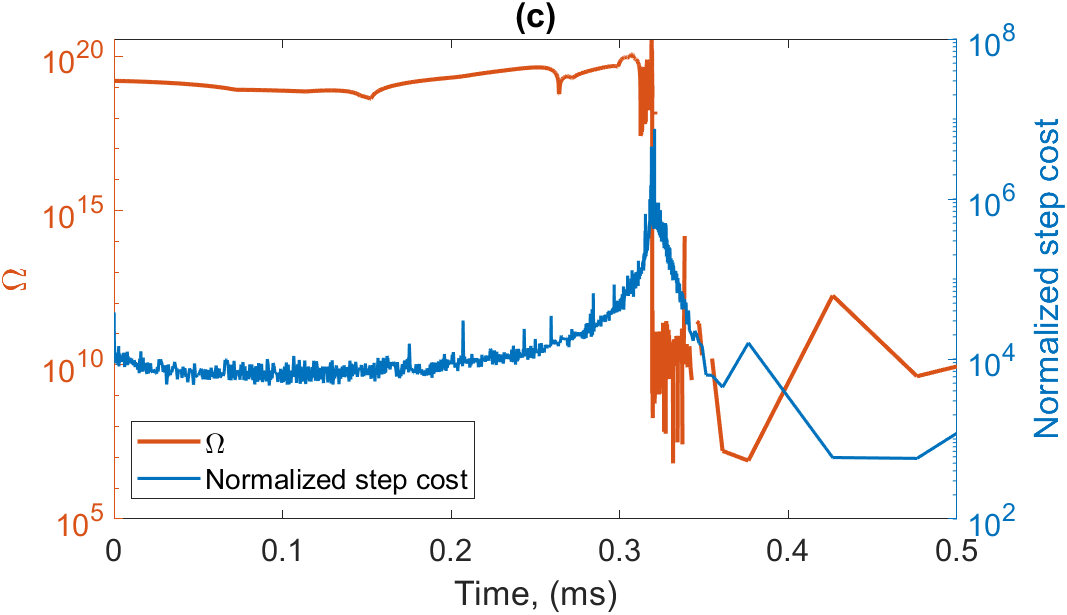}
			\end{subfigure}
			\begin{subfigure}[b]{0.45\textwidth}
				\centering
				\includegraphics[width=\textwidth]{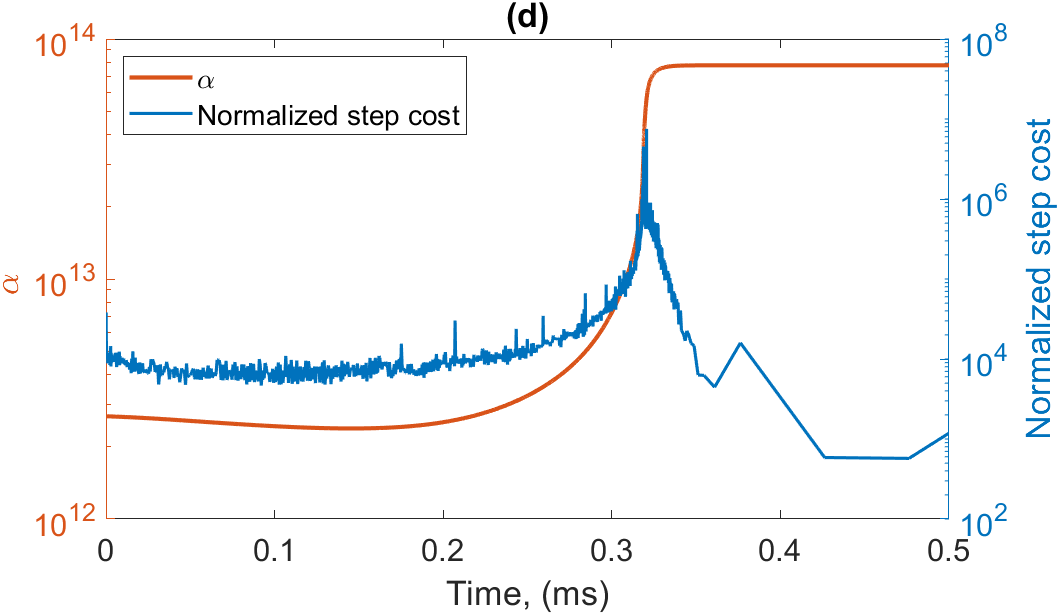}
			\end{subfigure}
			\begin{subfigure}[b]{0.45\textwidth}
				\includegraphics[width=\textwidth]{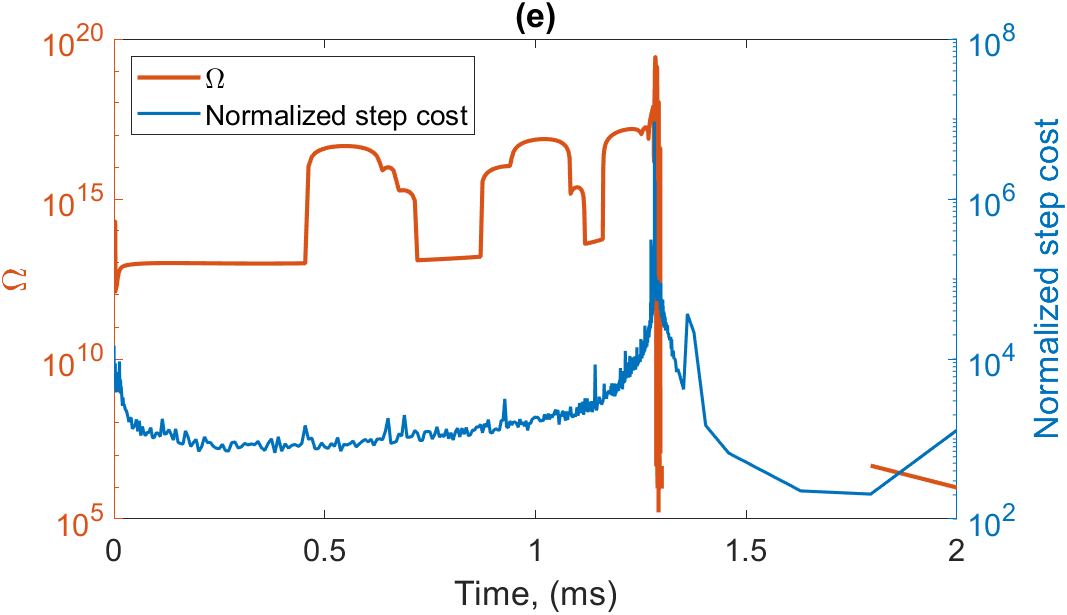} 
			\end{subfigure}
			\begin{subfigure}[b]{0.45\textwidth}
			\includegraphics[width=\textwidth]{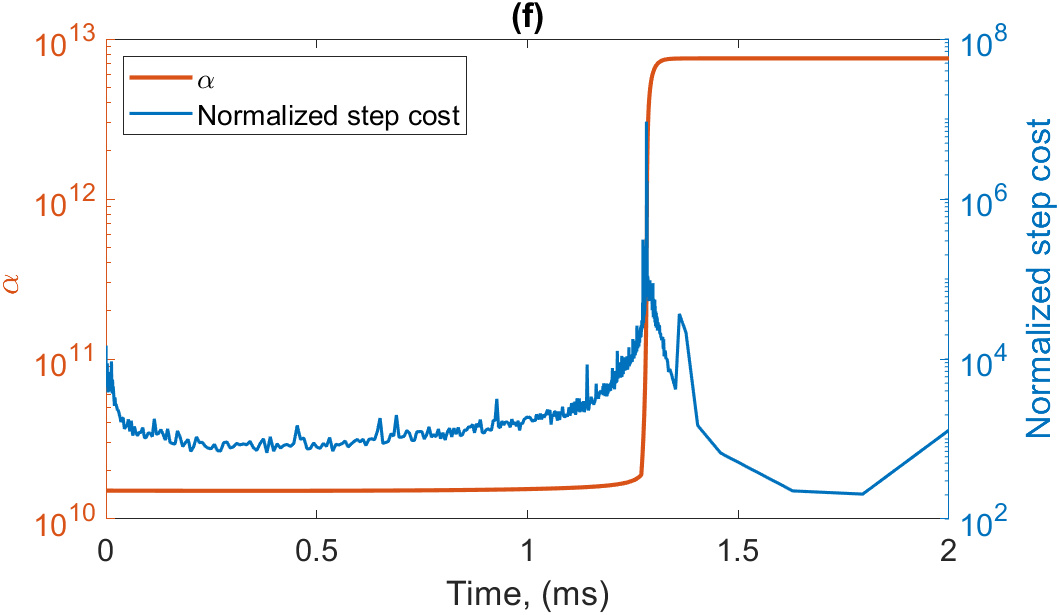} 
			\end{subfigure}
			
			\caption{  Plots vizualizing $\Omega$, which measures the area of the spectrum, and $\alpha$, measuring the real spread of the specturm, versus normalized step cost (the time spent integrating a step divided by the time step used).  Plots (a) and (b) plot information for the GRI mechanism. Plots (c) and (d) show the results for $\emph{n}$-dodecane, while  (e) and (f) demonstrates $\emph{n}$-butane data. $\Omega$ and $\alpha$ scales are set on the left axes, while the normalized step costs scales are on the right axes of the plots. }
			\label{SpecCost}
		\end{figure}
		For each experiment, we computed the eigenvalues of the exact Jacobian  provided by TCHEM using MATLAB's  {\tt eig} function.  The tightest set of tolerances from table (\ref{ExpTol}) were chosen, and the three experiments were carried out using the EPI3V time integration scheme, storing the Jacobian calculated by TCHEM at each step.  In order to provide a measure of the size of the spectrum of the Jacobian, we define the area $\Omega$ of the smallest rectangle with sides aligned with the coordinate axes such that it encloses all the eigenvalues in the complex plane, i.e. if $\lambda_j=a_j + \mathbf{i} \hspace{.1cm}b_j $ ($j = 1,...,N$) are the eigenvalues, the sides of the rectangle are $\alpha = \max_j(a_j)-\min_j(a_j)$ and $\beta = \max_j(b_j)-\min_j(b_j)$ so that $\Omega = \alpha \beta$.   Figure \ref{SpecCost} shows evolution of $\alpha$ and $\Omega$ during ignition, along with a normalized step cost, defined as the CPU time spent computing a step divided by the time interval stepped. 

	It is apparent in figure \ref{SpecCost} that the $\emph{n}$-dodecane mechanism has both the largest real spread $\alpha$ and $\Omega$ in all cases. 	 The KIOPS algorithm is based on projections onto the Krylov subspace and the estimation of exponentials of approximate eigenvalues.  If the problem's spectrum contains large positive real eigenvalues with a large $\Omega$, computing exponentials of these augmented systems is problematic; the adaptive time stepping procedure in KIOPS will reduce the time step size significantly to accommodate the user designated tolerance. However, in the case of the $\emph{n}$-dodecane mechanism, this time step reduction penalizes performance of EPIRK methods compared to the implicit scheme implemented in CVODE.  This increased cost is also reflected in the normalized CPU time in figure \ref{SpecCost}.

\section{Conclusion}\label{sec::conclusion}
  In our work, we investigated the performance of the novel EPI3V variable time stepping exponential integrator for the simulation of chemically reactive  and spatially homogeneous systems, i.e. chemical reactors.  We compared the performance of our EPI3V method to that of CVODE, which uses a modified Newton solver and Krylov-projection-type iterative method. Numerical expleriments were conducted for three chemical kinetics mechanisms of increasing complexity. We found that the exponential method performed favorably for certain mechanisms, but not for others. Comparable CPU time and accuracy were observed for both the GRI3.0 and $\emph{n}$-butane mechanism. However, for the $\emph{n}$-dodecane mechanism the CPU time for the EPI3V method required to obtain similar errors to the CVODE was an order of magnitude higher than that for CVODE at the same value of the error norm.

	We found the performance degradation of EPI3V method for $\emph{n}$-dodecane stems from a combination of a wide spectrum of the Jacobian and the presence of very large positive real eigenvalues. Because KIOPS is based on Krylov-iteration approximation of matrix exponential, its performance degrades in the presence of such spectra.  This finding points to a promising research direction to explore alternatives to Krylov-based algorithms for estimating products of matrix exponentials with vectors.  In the future we plan to investigate whether contour integration and quadrature-based methods will yield better performance \cite{Integral}. We are also currently extending our study to combustion problems that include transport, in particular to modeling flame front propagation. 
\section{Acknowledgement}
This work was done within the Exascale Catalytic Chemistry (ECC) Project, which is supported by the U.S. Department of Energy, Office of Science, Basic Energy Sciences, Chemical Sciences, Geosciences and Biosciences Division, as part of the Computational Chemistry Sciences Program. Sandia National Laboratories is a multimission laboratory managed and operated by National Technology $\&$ Engineering Solutions of Sandia, LLC, a wholly owned subsidiary of Honeywell International Inc., for the U.S. Department of Energy’s National Nuclear Security Administration under contract DE-NA0003525. This paper describes objective technical results and analysis. Any subjective views or opinions that might be expressed in the paper do not necessarily represent the views of the U.S. Department of Energy or the United States Government.

Additionally, this work was in part supported by the National Science Foundation under grants: 1840265, DMS-2012875.
\pagebreak

\bibliographystyle{IEEEtran}



\end{document}